\documentclass[12pt,a4paper,oneside]{article}
\usepackage{amsmath,amstext,amssymb,amsopn,amsthm}

\theoremstyle{plain}
\newtheorem{thr}{Theorem}

\newtheorem{lem}[thr]{Lemma}
\newtheorem{prop}[thr]{Proposition}

\theoremstyle{definition}
\newtheorem{defn}[thr]{Definition}

\theoremstyle{remark}
\newtheorem*{rem*}{Remark}
\newtheorem{rem}[thr]{Remark}

\newcommand{\CK}{the Chapman-Kolmogorov equation }
\newcommand{\Rd}{\mathbb{R}^d}

\newcommand{\Z}{\int^{\infty}_{0}}
\renewcommand{\leq}{\leqslant}
\renewcommand{\geq}{\geqslant}

\DeclareMathOperator{\diam}{diam}

\title{Intrinsic Ultracontractivity for L\'{e}vy Processes}
\author{Tomasz Grzywny\\
{\footnotesize Institute of Mathematics and Computer Science} \\
{\footnotesize Wroc\l{}aw University of Technology} \\
{\footnotesize Wyb{.} Wyspia\'nskiego 27, 50--370 Wroc\l{}aw, Poland}\\
{\footnotesize email: tomasz.grzywny@pwr.wroc.pl}} \date{}

\begin{document}
\maketitle
\begin{abstract}
We prove the intrinsic ultracontractivity for the semigroup
generated by a large class of  symmetric L\'{e}vy processes such
that the L\'{e}vy measure satisfies some conditions in the
neighborhood of $0$, killed on exiting a bounded and connected
Lipschitz domain.
\end{abstract}

\section{Introduction}
Intrinsic ultracontractivity  has been studied extensively in
recent year in the case of the symmetric diffusions (see e.g.
\cite{DS}, \cite{Banuelos}) and the symmetric $\alpha$-stable
process (see e.g. \cite{CS1}, \cite{K}). The concept of the
intrinsic ultracontractivity for non-symmetric semigroups was
introduced in \cite{KS}.

If the L\'{e}vy measure of symmetric L\'{e}vy processes $X_t$ is
"uniformly separate" from $0$ (see (\ref{Levy:this})) on truncated
cone with vertex in the neighborhood of $0$, we prove the
intrinsic ultracontractivity for semigroup generated by the killed
process on exiting a bounded and connected Lipschitz domain
(Theorem \ref{IUThr:this}). In the case if the Lebesgue measure is
absolutely continuous with respect to the L\'{e}vy measure then we
show that the semigroup is intrinsic ultracontractive for any
bounded open set (Remark \ref{remarkTh}).

The paper is organized in the following way. In Section 2 we
recall some definitions and prove facts about continuous and
strictly positivity of a transition density of process killed on
exiting a bounded open set. In Section 3 we prove the intrinsic
ultracontractivity.

\section{Preliminaries}
In  $\Rd$, $d\geq 1$, we consider a symmetric L\'{e}vy processes
$X_t$. By $\nu$ we denote its (nonzero) L\'{e}vy measure and by
$p(t,x,y)=p(t,x-y)$ the transition densities of $X_t$, which are
assumed to be continuous for every $t>0$ and defined for every $x,y\in \Rd
$. In addition  we assume
that there exists a constant $c(\delta)$ such that $p(t,x)\leq c$
for $t>0$ and $|x|\geq \delta$.

We use the notation $C=C(\alpha,\beta,\gamma,\dots)$ to denote
that the constant $C$ depends on $\alpha,\beta,\gamma,\dots$.
Usually values of constants may change from line to line, but they
are always  strictly positive and finite. Sometimes we skip in
notation that constants depend on usual quantities (e.g.
$d,D$). Next, we give some definitions. We denote
$$\tau_D=\inf\{t>0: X_t\notin D\},$$
$$\eta_D=\inf\{t\geq0: X_t\notin D\}.$$

Let $D$ be a bounded connected nonempty open set. In order to
study the killed process on exiting of $D$ we construct its
transition densities by the classical formula
$$p_D(t,x,y)=p(t,x,y)-r_D(t,x,y),$$
where
$$r_D(t,x,y)=E^x[t>\tau_D;p(t-\tau_D,X_{\tau_D},y)].$$
 The arguments used for Brownian motion (see
eg. \cite{CZ}) will prevail in our case and one can easily show
that $p_D(t,x,y), \ t\ge 0,$ satisfy the Chapman-Kolmogorov
equation (semigroup property). Moreover the transition density
$p_D(t,x,y)$ is a symmetric function $(x,y)$ a.s.. With the above
assumptions of the transition densities of the (free) process one
can actually show that $p_D(t,x,\cdot)$ and $p_D(t,\cdot,x)$ can
be chosen as continuous functions on $D$. The semigroup given by
the process $X_t$ killed on exiting of $D$ we denote by $P^D_t$.
We set $G_D(x,y)=\Z p_D(t,x,y)dt$ and call the Green function for
$D$.

$P^D_t$ is a strongly continuous semigroup of contractions on
$L^2(D)$. Because $p_D(t,x,y)$ is symmetric a.e., we obtain that
the operator $P^D_t$ is selfadjoint. For $D$ bounded we get from
continuity of $p(t,\cdot)$ that
\[p_D(t,x,y)\leq p(t,x-y)\leq \sup_{x\in B(0,\diam(D))}p(t,x)=C_1(t,D).\]
Therefore $P^D_t$ is Hilbert-Schmidt operator, so it's also
compact. So, it's well-known that there exists an orthonormal
basis of real-valued eigenfunctions $\{\varphi_n\}^{\infty}_{n=0}$
with corresponding eigenvalues $\{e^{-\lambda_n
t}\}^{\infty}_{n=0}$ satisfying
$0<\lambda_0<\lambda_1\leq\lambda_2\leq\ldots$, where all $\varphi_n$
are continuous.

We have that $p_D(t,x,\cdot)\in L^2(D)$ so we can represent this
function as
$$p_D(t,x,\cdot)=\sum^{\infty}_{n=0}<p_D(t,x,\cdot),\varphi_n>\varphi_n.$$
But $<p_D(t,x,\cdot),\varphi_n>=P^D_t\varphi_n(x)=e^{-\lambda_n
t}\varphi_n(x)$, so
$$p_D(t,x,y)=\sum^{\infty}_{n=0}e^{-\lambda_n
t}\varphi_n(x)\varphi_n(y).$$ Now, let us observe that the above
series are uniformly convergent, it follows from $|\varphi_n|\leq
e^{\lambda_n t/3}C_1(t/3,D)$ and
\[\sum^{\infty}_{n=0}e^{-\lambda_n
t/3}=\int_D p_D(t/3,x,x)dx\leq C_1(t/3,D)|D|.\] Hence, we get that
$p_D(t,\cdot,\cdot)\in C(D\times D)$. Therefore
$p_D(t,x,y)=p_D(t,y,x)$  for any $t>0$ and $x,y\in D$.

Next, we show that $p_D(t,\cdot,\cdot)$ is strictly positive on
$D\times D$. First, let us observe that for any $x\in D$ we have
$p_D(t,x,x)>0$. Indeed, $$p_D(t,x,x)=\int_D
p_D(\frac{t}{2},x,y)p_D(\frac{t}{2},y,x)dy =
\int_Dp^2_D(t/2,x,y)dy\geq (P^x(\tau_D>t/2))^2/|D|>0.$$ Let
$K\subset D$ be a compact and connected set. By continuity of
$p_D(t,\cdot,\cdot)$ we obtain that for any $x\in K$ there is a
radius $r_x$ such that
$$p_D(t,x,y)>0\quad \text{for } x,y\in B(x,2 r_x).$$ Because $K$
is compact, there are $x_1,\ldots,x_k\in K$ such that $K\subset
\bigcup_{i=1}^k B(x_k,r_{x_k})$. Now, we use a fact that $K$ is
connected to get from \CK that $p_D(k t,x,y)>0$ for any $x,y\in
K$. Hence we have that $p_D(s,x,y)>0$ for $s\geq kt$ and $x,y\in
K$.  Therefore $G_D(x,y)>0$, first for $x,y\in K$ and next for any
$x,y\in D$. This give us that $p_D(t,x,y)$ is strictly positive on
$D$ for any $t>0$. So we obtain that $\varphi_0$ is strictly
positive on $D$ too.

\begin{lem}\label{gestzabity2:this}
For any $x\in D$ and  $t>0$ we have
$$p_D(t,x,y)\leq C(t,D) E^x\tau_D E^y\tau_D.$$
\end{lem}

\begin{proof}
By the Chapman-Kolmogorov equation we obtain for $t> 0$
$$p_D(t,x,y)=\int_D p_D(t/2,x,z)p_D(t/2,z,y)dz \leq C_1(t/2,D) P^x(\tau_D>t/2).$$
Applying again the Chapman-Kolmogorov equation together with the
above inequality we get \begin{eqnarray*}p_D(t,x,y)&\leq& C_1
P^x(\tau_D>t/4)\int_D
p_D(t/2,z,y)dz\\&=&C_1P^x(\tau_D>t/4)P^y(\tau_D>t/2).\end{eqnarray*}
The application of Chebyshev's inequality completes the proof.
\end{proof}

\begin{defn} The semigroup $\{P^D_t\}$ is said to be {\em intrinsic
ultracontractive} if, for any $t>0$, there exists  a constant
$c_t$ such that
$$p_D(t,x,y)\leq c_t \varphi_0(x)\varphi_0(y),\quad x,y\in D.$$
\end{defn}

\begin{prop}\label{IUchar}
Let $D$ be a bounded connected nonempty open set. Then $\{P^D_t\}$
is intrinsic ultracontractive if and only if there is a constant
$C$ such that $E^x\tau_D\leq C \varphi_0(x)$.
\end{prop}
\begin{proof}
Suppose that $\{P^D_t\}$ is intrinsic ultracontractive that is
$$p_D(t,x,y)\leq c_t \varphi_0(x)\varphi_0(y).$$
Because $p_D(t,\cdot,\cdot)$ and $\varphi_0(\cdot)$ are continuous
and strictly positive, we have (see Theorem 3.2 in \cite{DS}) that there is
$\widetilde{c}_t$ such that
$$\widetilde{c}_t \varphi_0(x)\varphi_0(y)\leq  p_D(t,x,y).$$
If we integrate the above inequality with respect to $dt$ we get
$$C \varphi_0(x)\varphi_0(y)\leq  G_D(x,y).$$
And by integrating with respect to $dy$
$$\widetilde{C} \varphi(x)\leq  E^x\tau_D.$$

Now, suppose that $E^x\tau_D\leq C \varphi_0(x)$. From Lemma
\ref{gestzabity2:this} we have
$$p_D(t,x,y)\leq C_t E^x\tau_D E^y\tau_D,$$
what ends the proof.
\end{proof}

\section{Main results}
We prove intrinsic ultracontractivity for the semigroup
$P^{D}_{t}$ generated by the symmetric L\'{e}vy process, whose a
L\'{e}vy measure satisfies
\begin{equation}\label{Levy:this} \forall_{r>0,\gamma\in(0,\pi)}\exists_{\rho>0}
\inf_{|y|=\rho;\text{ }\Gamma_\gamma(y)}\nu(\Gamma_\gamma(y)\cap
B(0,r))>0,\end{equation} where $\Gamma_\gamma(y)$ is a right
circular cone of angle $\gamma$ at the vertex in $y$.

Notation and the proof of following theorem is similar as in paper
\cite{K}. We assume that $D$ is a bounded and connected Lipschitz
domain. That is there exist $\gamma_0$ and  $R_0>0$ and a cone
$\Gamma_{\gamma_0}=\{(y,x):0<x,\, y\in \mathbb{R}^{d-1},\,
\gamma_0 |y|< x\}$ such that for every $Q\in
\partial D$, there is a cone $\Gamma_{\gamma_0}(Q)$ with vertex $Q$, isometric
with $\Gamma_{\gamma_0}$ and satisfying $\Gamma_{\gamma_0}(Q)\cap
B(Q,R_0)\subset D$. Denote $U(\sigma)=\{x\in D:
\delta_D(x)<\sigma\}$, where
$\sigma\leq\frac{R_0}{4\sqrt{1+\gamma_0^2}}$. Then for any $x\in
U(\sigma)$ there are a point $y$ and a cone $\Gamma_{\gamma_0}(y)$
such that $|y-x|< \sigma(1+\sqrt{1+\gamma^2_0})\leq \frac{R_0}{2}$
and $\Gamma_{\gamma_0}(y)\cap B(x,R_0/2)\subset D\cap
U(\sigma)^c$.

We fix $x_0\in D$ and let $r>0$ be such that $\overline{B(x_0,2r)}
\subset D$. Denote $K=\overline{B(x_0,r)}$, $L=B(x_0,2r)$,
$M=D\backslash K$ and $N=D\backslash L$. We deal that $r\leq
\rho_0$. Define stopping time $S_n$ and $T_n$
\begin{eqnarray*}
S_1&=&0,\\
T_n&=&S_n+\eta_M\circ\theta_{S_n},\\
S_n&=&T_{n-1}+\eta_L\circ\theta_{T_{n-1}}.
\end{eqnarray*}

Now, we prove the following lemma.
\begin{lem}
There exists a constant $c=c(D,x_0)$ such that
\[P^x(X(\eta_M)\in K)\geq c E^x \eta_M\]
for all $x\in \Rd$.
\end{lem}

\begin{proof}
From (\ref{Levy:this}) we get existing a constant $\sigma_0\leq r$
such that
\[\inf_{|y|=\sigma_0(1+\sqrt{1+\gamma^2_0});\text{
}\Gamma_{\gamma_0}(y)}\nu(\Gamma_{\gamma_0}(y)\cap
B(0,R_0/2))=C_1.\] Denote $W=\{x\in D:\delta_D(x)\geq
\sigma_0/2\}\backslash B(x_0,r)$.

First, we prove that for $x\in W$, we have
\begin{equation}\label{lemma2_1:this}
P^x(X(\tau_M)\in K)\geq c_1,
\end{equation}
for some constant $c=c(r,D)$. Let $\rho_1$ be such that
\[\inf_{|y|=\rho_1;\text{
}\Gamma_{1}(y)}\nu(\Gamma_{1}(y)\cap B(0,r))=C_2>0.\] Denote
$J=D\backslash \overline{B(x_0,r-\rho_1/4)}$. Indeed, from the
Ikeda-Watanabe formula we have
\begin{eqnarray*}
P^x(X(\tau_M)\in K)&\geq& P^x\left(X(\tau_{J})\in
B(x_0,r-\rho_1/4)\right)\\
&\geq& P^x\left(X(\tau_{J})\in B(x_0,r-\rho_1/2)\right)\\
&=&\int_{J}G_J(x,y)\nu\left(B\left(x_0,r-\rho_1/2\right)-y\right)dy\\
&\geq&\int_{W}G_J(x,y)\nu\left(B\left(x_0,r-\rho_1/2\right)-y\right)dy.
\end{eqnarray*}
Because $p_J(t,\cdot,\cdot)$ is continuous and positive function
on $J\times J$ and $W\times W$ is compact subset of $J\times J$,
we get $\inf_{x,y\in W}p_J(t,x,y)>0$. So, $$\inf_{x,y\in
W}G_J(x,y)\geq\Z \inf_{x,y\in W}p_J(t,x,y)dt=c>0.$$ Besides, we
have $$\inf_{y\in
B(x_0,r+\rho_1/2)}\nu(\left(B\left(x_0,r-\rho_1/2\right)-y\right)\geq
\inf_{|y|=\rho_1;\text{ }\Gamma_{1}(y)}\nu(\Gamma_{1}(y)\cap
B(0,r))>0.$$ Therefore
\[P^x(X(\tau_M)\in K)\geq \varepsilon \int_{B(x_0,r+\rho_1/2)
\backslash
B(x_0,r)}\nu\left(B\left(x_0,r-\rho_1/2\right)-y\right)dy=C>0.\]
From (\ref{lemma2_1:this}) and the fact that $E^x\tau_M\leq
\widetilde{C}$ we obtain the claim of the lemma for $x\in W$.

Now, let $x\in D\backslash (W\cup K)=U(\sigma_0/2)$. Then from
Strong Markov Property we get
\begin{eqnarray*}
P^x(X(\tau_M)\in K) &=&
E^x(P^{X(\tau_{U(\sigma_0/2)})}(X(\tau_M)\in K))\geq c_2
E^x(E^{X(\tau_{U(\sigma_0/2)})}\tau_M)\\&=&c_2(E^x\tau_M-E^x\tau_{U(\sigma_0/2)}).
\end{eqnarray*}
And from (\ref{lemma2_1:this}) we obtain
\begin{eqnarray*}
P^x(X(\tau_M)\in K) &=& E^x(X(\tau_{U(\sigma_0/2)})\in W\cup
K,P^{X(\tau_{U(\sigma_0/2)})}(X(\tau_M)\in K))\\&\geq& c_1
P^x(X(\tau_{U(\sigma_0/2)})\in W\cup K).
\end{eqnarray*}
But \begin{eqnarray*}P^x(X(\tau_{U(\sigma_0/2)})\in W\cup K)&\geq&
P^x(X(\tau_{U(\sigma_0/2)}) \in D\cap
U(\sigma^0)) \\
&=& \int_{U(\sigma_0/2)}G_{U(\sigma_0/2)}(x,y)\nu(D\cap
U(\sigma^0)-y)dy\\&\geq& C_1
\int_{U(\sigma_0/2)}G_{U(\sigma_0/2)}(x,y)dy = C_1 E^x
\tau_{U(\sigma_0/2)}
\end{eqnarray*} Hence
\begin{eqnarray*}P^x(X(\tau_M)\in
K)&=&(\frac{1}{2}+\frac{1}{2})P^x(X(\tau_M)\in K)\\ &\geq&
\frac{c_2}{2}(E^x\tau_M-E^x\tau_{D\backslash (W\cup K)})+
\frac{C_1}{2}(E^x\tau_{D\backslash (W\cup K)})\\
&\geq& \frac{c_2\wedge C_1}{2} E^x\tau_M.\end{eqnarray*} For $x\in
D\backslash$ we have $E^x\tau_M=E^x\eta_M$, and the claim of the
lemma for $x\in D^c\cup K$ of course is obvious, so it ends the
proof.
\end{proof}

\begin{lem}
For all $x\in\Rd$ there exists a random variable $Z$ such that for
all $n\geq Z$ we have $T_n=\eta_D$ almost surely $P^x$.
\end{lem}

\begin{proof}
We will show that there exists a constant $\beta<1$ such that
$P^x(T_n<\eta_D)\leq \beta^{n}$ for all $x\in\Rd$ and $n\geq 1$.

Let $R=B(x_0,\diam(D))\backslash K$ and $\varepsilon=\inf_{x,y\in
B(x_0,\diam(D)-\rho_1/2)\backslash B(x_0,2r)}G_R(x,y)$ then from
the Ikeda-Watanabe formula for $x\in N$ we get
\begin{eqnarray*}
P^x(X(\eta_M)\in D^c)&\geq& P^x(X(\eta_R)\in B^c(x_0,\diam(D)))\\
&\geq&\int_{B(x_0,\diam(D)-\rho_1/2)\backslash
B(x_0,2r)}G_{R}(x,y)\nu(B^c(x_0-y,\diam(D)))dy\\&\geq& \varepsilon
\int_{B(0,\diam(D)-\rho_1/2)\backslash
B(0,\diam(D)-\rho_1)}\nu(B^c(y,\diam(D)))dy\geq \varepsilon C_2
c=1-\beta.
\end{eqnarray*}
Consequently, for any $x\in\Rd$ and $n\geq 1$ we get
$$
P^x(T_n<\eta_D, T_{n+1}=\eta_D)=P^x(T_n<\eta_D, S_{n+1}=\eta_D)+
P^x(T_n<\eta_D, S_{n+1}<\eta_D, X(T_{n+1})\in D^c)$$
\begin{eqnarray*}
&=& P^x(T_n<\eta_D, S_{n+1}=\eta_D)+
P^x(T_n<\eta_D, X(S_{n+1})\in N , X(\eta_M)\circ\theta_{S_{n+1}}\in D^c)\\
&=& P^x(T_n<\eta_D, S_{n+1}=\eta_D)+
E^x(T_n<\eta_D, X(S_{n+1})\in N , P^{X(S_{n+1})}(X(\eta_M)\in D^c))\\
&\geq& (1-\beta)P^x(T_n<\eta_D, S_{n+1}=\eta_D)+
(1-\beta)P^x(T_n<\eta_D,
S_{n+1}<\eta_D)\\&=&(1-\beta)P^x(T_n<\eta_D).
\end{eqnarray*}
Hence, we obtain
\begin{eqnarray*}
P^x(T_{n+1}<\eta_D)&=&P^x(T_n<\eta_D)-P^x(T_n<\eta_D,T_{n+1}=\eta_D)\\
&\leq& P^x(T_n<\eta_D) - (1-\beta)P^x(T_n<\eta_D) = \beta
P^x(T_n<\eta_D).
\end{eqnarray*}
Applying the Borel-Cantelli Lemma ends the proof of lemma.
\end{proof}

The above lemma allow us to prove similarly as Theorem 8 in
\cite{K} the following proposition.

\begin{prop}\label{propSetC}Let $C$ be an nonempty open subset of $D$. Then there
is $c$ such that $$ E^x\int^{\tau_D}_0\textbf{1}_C(X_t)dt\geq c
E^x\tau_D.$$
\end{prop}

\begin{thr}
There exists a constant $C$ such that
$$E^x\tau_D\leq C \varphi_0(x),$$
for all $x\in D$.
\end{thr}
\begin{proof}
We have, for all $t>0$,
$$e^{-\lambda_0 t}\varphi_0(x)=\int_D p_D(t,x,y)\varphi_0(y)dy.$$
By integration this with respect $dt$ we get
$$\varphi_0(x)=\lambda_0 \int_D G_D(x,y)\varphi_0(y)dy.$$
Because $\varphi_0$ is continuous and positive, we obtain that
there is a constant $\varepsilon>0$ such that a set $$C=\{x:
\varphi_0(x)>\varepsilon\}$$ is nonempty. By Proposition
\ref{propSetC} we have
\begin{eqnarray*}
E^x\tau_D&\leq& c^{-1}\int_CG_D(x,y)dy \leq
(c\varepsilon)^{-1}\int_C G_D(x,y)\varphi_0(y)dy\\ &\leq&
(c\varepsilon)^{-1}\int_D G_D(x,y)\varphi_0(y)dy =
(c\varepsilon\lambda_0)^{-1}\varphi_0(x).
\end{eqnarray*}
\end{proof}
Applying Lemma \ref{IUchar} give us the theorem below.
\begin{thr}\label{IUThr:this}
Let $D$ be an bounded and connected Lipschitz domain. If the
L\'{e}vy measure of symmetric L\'{e}vy process $X_t$ satisfies
(\ref{Levy:this}), then the semigroup $\{P^D_t\}$ is intrinsic
ultracontractive.
\end{thr}

\begin{rem}\label{remarkTh}
Suppose that the symmetric L\'{e}vy process $X_t$ has the L\'{e}vy
measure such that the Lebesgue measure is absolutely continuous
with respect to it. Then the semigroup $P^D_t$ is intrinsic
ultracontractive for any bounded open set.
\end{rem}
\begin{proof}
Proof of this remark is the same as the proof of Theorem 1 in
\cite{K}.
\end{proof}

\end{document}